\documentclass[letterpaper,12pt]{article}
\usepackage{times}
\usepackage{empheq}
\usepackage{amsmath, amssymb}
\usepackage{amsthm}
\usepackage[T1]{fontenc}
\usepackage[english]{babel}
\usepackage[utf8]{inputenc}
\usepackage{hyperref}
\usepackage{mathrsfs}
\usepackage{mathabx}
\usepackage{color}
\usepackage{anysize}
\usepackage{apacite}
\usepackage{natbib}
\marginsize{2.5cm}{2.5cm}{2.5cm}{2.5cm}
\usepackage[toc,page]{appendix}
\usepackage{url}
\usepackage{mathtools}
\usepackage[font={small,sc}]{caption}
\usepackage{ragged2e}

\newtheorem{Theo}{Theorem}[section]

\newtheorem{Prop}{Proposition}[section]
\newtheorem{Corol}{Corollary}[section]

\newtheorem{Rem}{Remark}[section]

\AtBeginDocument{}

\DeclareMathOperator\Span{span}

\newcommand{\Rd}{\mathbb{R}^{d}}

\newcommand{\RdxRd}{\Rd\times\Rd}

\newcommand{\SchwartzRd}{\mathscr{S}(\Rd)}

\newcommand{\SchwartzRdxRd}{\mathscr{S}(\RdxRd)}
\newcommand{\TemperedRd}{\mathscr{S}'(\Rd)}

\newcommand{\TemperedRdxRd}{\mathscr{S}'(\RdxRd)}

\newcommand{\Cov}{\mathbb{C}ov}
\newcommand{\Var}{\mathbb{V}ar}
\newcommand{\OmegaAP}{(\Omega , \mathcal{A} , \mathbb{P})}
\newcommand{\LtwoOmega}{L^{2}(\Omega , \mathcal{A} , \mathbb{P})}

\newcommand{\DRd}{\mathscr{D}(\Rd)}

\newcommand{\DistributionsRd}{\mathscr{D}'(\Rd)}

\newcommand{\LtwoRd}{L^{2}(\Rd)}
\newcommand{\LtwoRdmu}{L^{2}(\Rd,\mu)}

\begin{document}

\begin{Large}
\begin{center}
\textbf{Generalized Stochastic Processes: Linear Relations to White Noise and Orthogonal Representations}
\end{center}
\end{Large}
\begin{center}
\textsc{Carrizo Vergara, Ricardo}
\end{center}
\begin{center}
Université de Paris II Panthéon-Assas, Paris, France. \\
\textit{Ricardo.Carrizo-Vergara@u-paris2.fr} 
\end{center}
\vspace{-3ex}
\begin{center}
2021
\end{center}

\begin{center}
\textbf{Abstract}
\end{center}

\small
\justify

We present two linear relations between an arbitrary (real tempered second order) generalized stochastic process over $\Rd$ and White Noise processes over $\Rd$. The first is that any generalized stochastic process can be obtained as a linear transformation of a White Noise. The second indicates that, under dimensional compatibility conditions, a generalized stochastic process can be linearly transformed into a White Noise. The arguments rely on the regularity theorem for tempered distributions, which is used to obtain a mean-square continuous stochastic process which is then expressed in a Karhunen-Loève expansion with respect to a convenient Hilbert space. The first linear relation obtained allows also to conclude that any generalized stochastic process has an orthogonal representation as a series expansion of deterministic tempered distributions weighted by uncorrelated random variables with summable variances. This representation is then used to conclude the second linear relation.
\bigskip


\normalsize

\section{Introduction}
\label{Sec:Introduction}

Spatial stochastic processes are a fundamental tool in both practice and theory of Stochastic Analysis. An elegant generalization of regular (mean-square continuous, for instance) stochastic processes is the one of Generalized Stochastic Process (GeSP), that is a random distribution. Many authors have used and studied such concept \citep{gelfand1964generalized,ito1954stationary,matheron1965variables,rozanov1982markov}. One important advantage is the possibility of defining, without mayor difficulties, derivatives and Fourier transforms for large classes of processes. A Second order GeSP over $\Rd$ has its statistical properties partially determined (sometimes completely) by its covariance structure, which may described either by a sesquilinear form over the space of test-functions (the covariance Kernel) or by a covariance distribution. A particularly important case is the one of White Noise, whose covariance Kernel is given by the inner product of $\LtwoRd$. This has many consequences which make White Noise the cornerstone of many theoretical developments and practical methodologies. To list a few, White Noise is an orthogonal mean-square random measure with translation invariant variance, so it is easy to simulate; the Fourier transform of a White Noise has also the structure of a White Noise, which has theoretical and practical consequences in Signal Analysis; and White Noise is the distributional derivative of Brownian motion in $d=1$ (and the $\frac{\partial^{d}}{\partial x_{1} ... \partial x_{d}}$-derivative of Brownian sheet in dimension $d$), so basically all modern Stochastic Analysis can be grounded on it. We refer to 
\citep{walsh1986introduction,holden2009stochastic,hida2013white} as important examples of theoretical developments around White Noise.

In this work we prove that any real second order GeSP can be expressed as a linear transformation of a White Noise over $\Rd$. This result may be seen as a generalization to the infinite (countable) dimensional case of the Cholesky method of simulating random vectors, or of Principal Components Analysis, applied to the case of random distributions. The argument, while requiring some technicality proper of the use of tempered distributions and topological vector spaces in general, is actually quite simple. If $Z$ is a real GeSP it can be concluded from an application of the regularity theorem for tempered distributions that $Z$ is the derivative (we use a fractional Laplacian for convenience) of high-enough order of a mean-square continuous stochastic process with polynomially bounded covariance, say $Y$. This $Y$ can be then conceived as a process acting linearly and mean-square continuously over the Hilbert space $L^{2}(\Rd,\mu)$, being $\mu$ a convenient finite measure over $\Rd$. A Karhunen-Loève expansion can be then obtained for $Y$, in terms of the eigenfunctions and eigenvalues of its  associated integral covariance operator. We re-obtain then $Z$ by deriving $Y$ enough times. It follows that an orthogonal representation for $Z$ can be given:
\begin{equation}
\label{Eq:OrthogonalZIntro}
Z = \sum_{n \in \mathbb{N}}Z_{n}T_{n},
\end{equation}
where $(T_{n})_{n \in \mathbb{N}}$ is a sequence of linearly independent tempered distributions and $(Z_{n})_{n \in \mathbb{N}}$ is a sequence of uncorrelated random variables whose variances form a convergent series. Similar orthogonal representations are known for processes having values in Hilbert spaces or their duals and for particular cases of GeSPs \citep{meidan1979reproducing,angulo1997orthogonal,angulo2006wavelet,medina2005convergence}; here we propose a construction of such representation applicable to any (tempered) random distribution, and where the distributions $T_{n}$ are obtained explicitly, without concerning about Hilbert subspaces. Since a White Noise over $\Rd$ can be expressed as
\begin{equation}
\label{Eq:OrthogonalWIntro}
W = \sum_{n \in \mathbb{N}}\epsilon_{n}g_{n},
\end{equation}
with $(\epsilon_{n})_{n \in \mathbb{N}}$ uncorrelated with unit variances and $(g_{n})_{n \in \mathbb{N}}$ an orthonormal basis of $\LtwoRd$, we can construct a White Noise using the variables $Z_{n}$ normalized and some basis of $\LtwoRd$. Then, we re-obtain $Z$ by applying a convenient operator which maps the basis of $\LtwoRd$ to the distributions $T_{n}$ (up to a multiplicative constant). Hence, 
\begin{equation}
\label{Eq:Z=LWIntro}
Z = \mathcal{L}W
\end{equation}
for some White Noise $W$ and for some linear operator $\mathcal{L}$ which \textit{can be applied to a White Noise}. On the other hand, if there is an infinity quantity of the random variables $(Z_{n})_{n \in \mathbb{N}}$ with non-null variance, by mapping conveniently the distributions $T_{n}$ to a basis of $\LtwoRd$, we obtain an equation of the form 
\begin{equation}
\label{Eq:LZ=WIntro}
\mathcal{L}Z = W,
\end{equation} 
for some (other) White Noise $W$ and some (other) linear operator $\mathcal{L}$ \textit{applicable to $Z$}. This answers, in a particularly abstract way, a question which has been made in the practical framework of the SPDE-Approach in Spatial Statistics: given a random field  over $\Rd$ with a given covariance, when can we find a White Noise (additively) driven SPDE that the random field satisfies? See \citet{lindgren2011explicit} for the origins of the SPDE-Approach and the conclusive part of \citet{carrizov2018development} for the mentioned question.

Relations \eqref{Eq:Z=LWIntro} and \eqref{Eq:LZ=WIntro} are then the linear relations of $Z$ to White Noise which are discussed in this work, and the arguments to those are essentially completely given in this introduction. The main of the corpus of this article consists then in the rigorous details of every argument used, which includes the precise definitions of a GeSP and a White Noise over $\Rd$, when and how to apply continuous linear operators over them and in which sense those are continuous. It also includes an explicit construction of the distributions $(T_{n})_{n \in \mathbb{N}}$ and, derived from it, the explicit construction of a Hilbert space associated to $Z$ which is used to properly obtain the result \eqref{Eq:LZ=WIntro}, specifying the domain of application of the operator $\mathcal{L}$ in such case.

After specifying some notations and general settings, in Section \ref{Sec:KL} we present the version of the Karhunen-Loève Theorem which is useful for us, together with a proof. In Section \ref{Sec:GeSPannWN} we give an overview on GeSPs, and we specify how and when can we apply linear operators on them. In Section \ref{Sec:RepresentationZ=LW} we present the linear transformation from White Noise and the orthogonal representation result. In Section \ref{Sec:LZ=W} we transform, when we can, a GeSP into a White Noise. We finish in Section \ref{Sec:Conclusion} with some concluding remarks.

\subsection*{General settings and notations}

Every complex vector space considered in this work is supposed to be endowed with a complex conjugation with respect to which it is closed, the complex conjugate of an element $e$ being denoted $\overline{e}$. For a vector space $E$ we denote $E^{*}$ its algebraic dual and $E'$ its topological dual if $E$ is topological. If $T \in E^{*}$, its action over an element $e \in E$ is denoted by $\langle T , e \rangle$. Over $E'$, the topology of pointwise convergence is called the $E'$-weak topology (also called the weak-$*$ topology). If $E$ is a normed space we denote its norm by $\| \cdot \|_{E}$ and if it is Hilbert, its inner product is denoted $(\cdot , \cdot)_{E}$ and it is taken to be a sesquilinear, being antilinear in the second variable. In addition, if $e \in E$, we denote $e^{*} \in E'$ the linear functional such that $ \langle e^{*} , \cdot \rangle = ( \cdot , \overline{e} )_{E}$. In Appendix \ref{Sec:AppReminders} we give some reminders on locally convex topological vector spaces and topological tensor products.

$\SchwartzRd$ denotes the Schwartz space of rapidly decreasing complex smooth functions over $\Rd$. When endowed with its classical topology (see Appendix \ref{Sec:AppReminders} for a reminder) its dual space $\TemperedRd$ is the space of tempered distributions. $C_{PB}(\Rd)$ denotes the space of continuous and polynomially bounded complex functions over $\Rd$. If $f$ is an enough regular complex function over $\Rd$ (in $C_{PB}(\Rd)$ or in $L^{p}(\Rd)$ with $p \in \left[ 1 , \infty \right]$, for instance), we denote by $f^{*} \in \TemperedRd$ its associated tempered distribution, that is $\langle f^{*} , \varphi \rangle = \int_{\Rd} f(x)\varphi(x)dx$ for every $\varphi \in \SchwartzRd$.

We will work with random variables defined over a fixed probability space $\OmegaAP$. All the random objects are supposed to have zero mean. By complex \textit{stochastic process} we mean a family of complex random variables $(X(t))_{t \in T}$ indexed by a non-empty set $T$. The equality between two stochastic processes is always interpreted as one being a modification of the other. We do not focus on the properties of the trajectories of the processes, but mainly on its mean-square behaviour. The convergence of random objects considered in this work are always in a mean-square sense (we make it precise every time we consider it). We do not specify the laws of the random variables involved, but we assume that the used processes can be formally constructed (as it is the case for Gaussian processes, for instance, using Kolmogorov's Theorem).

\section{Karhunen-Loève Expansion}
\label{Sec:KL}

The typical use of the Karhunen-Loève expansion is done in the case of mean-square continuous stochastic processes over compact domains. In such a case the convergence of the expansion is in mean-square and uniformly in the indexation set \citep[Section 37.5]{loeve1978probability}. Other versions suppose only that the trajectories of the process belong to a particular Hilbert space \citep{red2009elements}. From our part, we will use a more general and weak version of the Karhunen-Loève Theorem applicable to a stochastic process indexed linearly and continuously by a separable Hilbert space, without making precise if the trajectories of the process are or not in the topological dual of such space. In such a case, a mean-square-weak-$*$ convergence is all we actually need.

\begin{Theo}[\textbf{Karhunen-Loève}]
\label{Theo:KL}
Let $E$ be a separable complex Hilbert space. Let $(X(e))_{e \in E}$ be a complex stochastic process indexed by $E$ satisfying the following properties:
\begin{enumerate}
\item \label{Itm:1KL} $X(e) \in \LtwoOmega$ for all $e \in E$.
\item \label{Itm:2KL} $X(\alpha e + \beta f ) \stackrel{a.s.}{=} \alpha X(e) + \beta X(f) $, for all $e,f \in E$ and for all $\alpha,\beta \in \mathbb{C}$.
\item \label{Itm:3KL} The covariance sesquilinear kernel $K_{X}(e,f) = \mathbb{E}(X(e)\overline{X(f)})$ is bounded, and the induced covariance operator $Q_{X} : E \to E $ is compact.
\end{enumerate}
Then, there exists a complete orthonormal basis of $E$, $(f_{n})_{n \in \mathbb{N}} \subset E $ such that
\begin{equation}
\label{Eq:X(e)KLexpansion}
X(e) \stackrel{a.s.}{=} \sum_{n \in \mathbb{N}} X(f_{n})( e , f_{n} )_{E}, \quad \forall e \in E,
\end{equation}
where the series is taken in $\LtwoOmega$ and the random variables $(X(f_{n}))_{n \in \mathbb{N}}$ are pairwise uncorrelated with converging-to-zero variances.
\end{Theo}

Equation \eqref{Eq:X(e)KLexpansion} can also be expressed as
\begin{equation}
\label{Eq:XKLexpansion}
X = \sum_{n \in \mathbb{N}}X(f_{n})\overline{f_{n}}^{*},
\end{equation}
where the convergence of the series is understood in a mean-square-$E'$-weak sense.  The explicit definition of the \textit{covariance operator} $Q_{X}$ is given in the proof.

\textbf{Proof: } From  conditions \ref{Itm:1KL}, \ref{Itm:2KL}, \ref{Itm:3KL}, the covariance kernel $K_{X} : E \times E \to \mathbb{C} $ given by $K_{X}(e,f) = \mathbb{E}( X(e)\overline{X(f)} )$ is a well-defined positive semidefinite and continuous sesquilinear form on $E$. There exists then $C > 0 $ such that
\begin{equation}
\label{Eq:Kbounded}
|K_{X}(e,f)| \leq C \| e\|_{E} \| f \|_{E}, \quad \forall e,f \in E.
\end{equation}
Hence, for every $f \in E$, the functional $K_{X}( \cdot , f ) $ is in $E'$. Its norm satisfies $\| K_{X}( \cdot , f ) \|_{E'} \leq C \| f \|_{E}$. By Riesz Representation Theorem, for every $f \in E$ there exists a unique $ Q_{X}(f) \in E $ identifying the linear functional $K_{X}(\cdot , f) $. The so-defined operation $Q_{X} : E \to E$ satisfies
\begin{equation}
\label{Eq:K(e,f)=(e,Qf)}
K_{X}(e,f) = (e , Q_{X}(f) )_{E}, \quad \forall e,f \in E.
\end{equation} 
Using the sesquilinearity of both the inner product and $K_{X}$, one proves easily that $Q_{X}$ is linear. One also has $\| Q_{X}(f) \|_{E}  = \| K_{X}(\cdot , f ) \|_{E'} \leq C \| f \|_{E}$, hence $Q_{X}$ is bounded. $Q_{X}$ is called the \textit{covariance operator} of $X$. Because of the relation \eqref{Eq:K(e,f)=(e,Qf)}, it is clear that $Q_{X}$ is positive semidefinite and self-adjoint.

By supposition \ref{Itm:3KL}, $Q_{X}$ is compact. Since it is also self-adjoint and $E$ is separable, the Hilbert-Schmidt Theorem \citep[Theorem VI.16]{reed1980methods} allows us to conclude that there exists an orthonormal basis of $E$, $(f_{n})_{n \in \mathbb{N}}$ consisting of eigenvectors of $Q_{X}$, for which the associated eigenvalues $(\lambda_{n})_{n \in \mathbb{N}}$ are non-negative ($Q_{X}$ is positive semidefinite) and converging to zero. It is then concluded that the sequence of random variables $( X(f_{n}) )_{n \in \mathbb{N}}$ satisfies
\begin{equation}
\label{Eq:X(fn)orthogonal}
\mathbb{E}\left(X(f_{n})\overline{X(f_{m})} \right) = \left( f_{n} , Q_{X}\left(f_{m}\right) \right)_{E} = \lambda_{n} \delta_{n,m},
\end{equation}
hence they are pairwise uncorrelated with converging-to-zero variances. 

Finally, since $ E( | X(e) |^{2} ) = \left( e , Q_{X}(e) \right)_{E} \leq C \| e \|_{E}^{2}$, one has that $X$ can be interpreted as a bounded linear operator $ X : E \to \LtwoOmega$. It follows that if we have a sequence $(e_{n})_{n} \subset E$ converging to $e \in E$, then $X(e_{n}) \to X(e) $ in $\LtwoOmega$. By using the expansion of $e$ in the basis $(f_{n})_{n \in \mathbb{N}}$, one obtains
\begin{equation}
\label{Eq:X(e)InBasisfn}
X(e) = X \left( \sum_{n \in \mathbb{N}} (e , f_{n})_{E}f_{n} \right) \stackrel{a.s.}{=} \sum_{n \in \mathbb{N}} X(f_{n})(e,f_{n})_{E},
\end{equation}
where the series is taken in the sense of $\LtwoOmega$. $\blacksquare$

\section{Generalized Stochastic Processes and White Noise}
\label{Sec:GeSPannWN}

\subsection{Generalized Stochastic Processes}
\label{Sec:GeSP}

A second order real (tempered) Generalized Stochastic Process over $\Rd$ (or Generalized Random Field, or Random Distribution, abbreviated GeSP from now on) is a real linear and continuous mapping $Z : \SchwartzRd \to \LtwoOmega$. By \textit{real}, we mean that when applied to real test-functions, the result is a real random variable. A GeSP can be conceived as a stochastic process indexed by the Schwartz space, $( Z(\varphi) )_{\varphi \in \SchwartzRd}$, having a linear and mean-square continuous behaviour over $\SchwartzRd$. The covariance Kernel of $Z$ is defined as the mapping $K_{Z} : \SchwartzRd \times \SchwartzRd \to \mathbb{C}$ given by $K_{Z}$ is a sesquilinear positive semidefinite form on $\SchwartzRd$. From Cauchy-Schwarz inequality and the continuity of $Z$ one concludes that $K_{Z}$ is separately continuous.  From Schwartz's Nuclear Theorem (see 
\citet[Theorem V.12]{reed1980methods} 
 or \citet[Theorem 51.6 and its Corollary]{treves1967topological}), there exists a unique real tempered distribution $C_{Z} \in \TemperedRdxRd$, called the \textit{covariance distribution}, such that
\begin{equation}
\label{Eq:CZdef}
\langle C_{Z} , \varphi\otimes\overline{\phi} \rangle = K_{Z}(\varphi , \phi) , \quad \forall\varphi,\phi \in \SchwartzRd,
\end{equation}
where $\varphi \otimes \phi $ is the tensor product function $ ( \varphi \otimes \phi )(x,y) = \varphi(x)\phi(y) $. $C_{Z}$ is positive semidefinite in the sense that $ \langle C_{Z} , \varphi \otimes \overline{\varphi} \rangle \geq 0$ for all $\varphi \in \SchwartzRd$.

When the covariance distribution $C_{Z}$ can be identified with a continuous function $f_{C_{Z}} : \RdxRd \to \mathbb{R}$, that is, when
\begin{equation}
\label{Eq:CZcontinuous}
\langle C_{Z} ,  \psi \rangle = \int_{\RdxRd} f_{C_{Z}}(x,y) \psi(x,y) d(x,y), \quad \forall \psi \in \SchwartzRdxRd,
\end{equation}
then the GeSP $Z$ can also be identified with a stochastic process indexed by $\Rd$, $(Z(x))_{x \in \Rd}$ which is mean-square continuous with covariance function $f_{C_{Z}}$. The random variable $Z(x)$ can be defined as a limit in mean-square of sequences of the form $(Z(\varphi_{n}))_{n \in \mathbb{N}}$ with $\varphi_{n} \to \delta_{x}$ in a particular sense (for example, using mollifiers). It can also be proven that
\begin{equation}
\label{Eq:Z(phi)Integral}
\int_{\Rd}Z(x)\varphi(x)dx \stackrel{a.s.}{=} Z(\varphi), \quad \forall \varphi \in \SchwartzRd,
\end{equation}
expression which allows to re-obtain the random variables of the GeSP $(Z(\varphi))_{\varphi \in \SchwartzRd}$ starting from its regular version $(Z(x))_{x \in \Rd}$.

\subsection{Operators over GeSPs}
\label{Sec:OperatorsGeSP}

If $\mathcal{L} : \SchwartzRd \to \SchwartzRd$ is a continuous and linear operator, its adjoint $\mathcal{L}^{*} : \TemperedRd \to \TemperedRd $ is weak-$*$ continuous and can be applied to any GeSP $Z$ through the formula
\begin{equation}
\label{Eq:L*Z}
\mathcal{L}^{*}Z (\varphi) := Z ( \mathcal{L}\varphi ), \quad \forall\varphi \in \SchwartzRd.
\end{equation}
The so-defined process $\mathcal{L}^{*}Z$ is also a GeSP with covariance distribution \citep[Section 3.4.2]{carrizov2018development}
\begin{equation}
\label{Eq:CovL*Z}
C_{\mathcal{L}^{*}Z} = ( \mathcal{L}^{*}\otimes\overline{\mathcal{L}^{*}} ) C_{Z}.
\end{equation}
This principle can be applied to define the application of derivatives and fractional differential operators over GeSPs without technical difficulties. This is the case for example of fractional Laplacian operators. We recall that for $\alpha \in \mathbb{R}$, the operator $(1 - \Delta)^{\alpha}$ over $\TemperedRd$ can be defined as $(1-\Delta)^{\alpha} := \mathscr{F}^{-1}(  (1+|\xi|^{2})^{\alpha}\mathscr{F}(\cdot)  ) $, where $\mathscr{F}$ denotes de Fourier Transform over $\Rd$. This operator is well-defined over $\TemperedRd$ and it is bijective with $[ (1 - \Delta)^{\alpha} ]^{-1} = (1 - \Delta)^{-\alpha}$.

We will use the following important result on tempered distributions \citep[Theorem V. 10]{reed1980methods}.
\begin{Theo}[\textbf{Regularity Theorem for tempered distributions}]
\label{Theo:RegularityTempered}
Let $T \in \TemperedRd$. Then, $T$ is a derivative of high-enough order of a function in $C_{PB}(\Rd)$. In particular, there exists $N \in \mathbb{N}$ such that the distribution $(1 - \Delta )^{-\frac{N}{2}}T$ can be identified with a function in $C_{PB}(\Rd)$.
\end{Theo}

From this theorem we can conclude the following result:
\begin{Prop}[\textbf{Regularity Theorem for GeSPs}]
\label{Prop:RegularityGeSP}
Let $Z$ be a real GeSP over $\Rd$. Then, there exists $N \in \mathbb{N}$ such that the GeSP $(1-\Delta)^{-\frac{N}{2}}Z$ can be identified with a mean-square continuous stochastic process with polynomially bounded covariance function.
\end{Prop}

\textbf{Proof:} Since by Nuclear Theorem the covariance distribution $C_{Z}$ is in $\TemperedRdxRd$, by a suitable adaptation of regularity Theorem \ref{Theo:RegularityTempered}  one can find a big-enough $N \in \mathbb{N}$ such that $ (1-\Delta_{x})^{-\frac{N}{2}}(1-\Delta_{y})^{-\frac{N}{2}}C_{Z}$ is in $C_{PB}(\RdxRd)$. By \eqref{Eq:CovL*Z}, this distribution coincides with the covariance distribution of $(1-\Delta)^{-\frac{N}{2}}Z$, hence this GeSP is identified with a mean-square continuous stochastic process with covariance function in $C_{PB}(\RdxRd)$. $\blacksquare$    

\subsection{Operators over subclasses of GeSPs}

The operators considered until now can be freely applied to any GeSP. However, if we consider a linear operator $\mathcal{L}$ defined only in a subspace $\mathscr{V}$ of $\TemperedRd$ then its application to an arbitrary GeSP is not immediate. Nevertheless, it could be applied to some GeSP which ``behaves as a random member of $\mathscr{V}$''. As an example, consider the Dirac operator $\delta_{x} : C_{PB}(\Rd) \to \mathbb{C}$ given by $\langle \delta_{x}, f \rangle = f(x) $. Then, $\delta_{x}$ cannot be immediately applied to an arbitrary GeSP (the punctual evaluation may not make sense for a random distribution), but it can be properly applied to every GeSP whose covariance is in $C_{PB}(\RdxRd)$. In such case, it is the regularity of the covariance distribution which allows us to apply the operator.

We will make precise this notion with the following general Proposition. We use the concept of projective topological tensor product of two locally convex topological vector spaces, see Appendix \ref{Sec:AppReminders}. The key issue is to ask the covariance to belong to a particular subspace of $\TemperedRdxRd$.

\begin{Prop}
\label{Prop:ExtenstionZtoU}
Let $\mathscr{U}(\Rd)$ be a complex Hausdorff locally convex topological vector space such that $\SchwartzRd $ is (identified with) a dense subspace of $\mathscr{U}(\Rd)$ with continuous embedding. Let $Z$ be a real GeSP over $\Rd$ such that its covariance $C_{Z}$ is in $( \mathscr{U}(\Rd) \otimes_{\pi} \mathscr{U}(\Rd) )'$. Then,
\begin{itemize}
\item $Z$ can be uniquely extended to a stochastic process which is a linear and continuous mapping from $\mathscr{U}(\Rd)$ to $\LtwoOmega$.
\item If $\mathcal{L} : \SchwartzRd \to \mathscr{U}(\Rd)$ is linear and continuous, then its adjoint $ \mathcal{L}^{*} : \mathscr{U}'(\Rd) \to \TemperedRd$ can be applied to $Z$ through the formula
\begin{equation}
\label{Eq:L*Zextended}
\mathcal{L}^{*}Z (\varphi) = Z( \mathcal{L}\varphi  ), \quad \forall \varphi \in \SchwartzRd,
\end{equation}
the result being a well-defined GeSP over $\Rd$.
\end{itemize}
\end{Prop}

The proof of this Proposition is given in Appendix \ref{Sec:AppReminders}. Let us give a Corollary for the case where $\mathscr{U}(\Rd)$ is a Hilbert space, which will be of particular importance and which can be directly obtained from the positive definiteness of the covariance distribution.

\begin{Corol}
\label{Corol:UHilbert}
Proposition \ref{Prop:ExtenstionZtoU} holds if $\mathscr{U}(\Rd)$ is a Hilbert space and if the covariance distribution of $Z$ satisfies $\langle  C_{Z},  \varphi \otimes \overline{\varphi} \rangle  \leq  \| \varphi \|^{2}_{\mathscr{U}(\Rd)}$ for every $\varphi \in \SchwartzRd$.
\end{Corol}

\subsection{White Noise}
\label{Sec:WN}

The most important example of GeSP is the White Noise over $\Rd$, noted $W$. A real White Noise over $\Rd$ is a real GeSP with covariance distribution given by the $L^{2}(\Rd)$ inner product:
\begin{equation}
\label{Eq:CovW}
 \langle C_{W} , \varphi \otimes \overline{\phi} \rangle = \left(\varphi , \phi \right)_{L^{2}(\Rd)} = \int_{\Rd}\varphi(x)\overline{\phi(x)}dx, \quad \forall \varphi, \phi \in \SchwartzRd.
\end{equation}
Following Corollary \ref{Corol:UHilbert}, the domain of indexation of $(W(\varphi))_{\varphi \in \SchwartzRd}$ can be extended continuously to $\varphi \in \LtwoRd$ using the density and the continuous embedding of $\SchwartzRd$ in $\LtwoRd$. It follows that if $(e_{n})_{n \in \mathbb{N}}$ is an orthonormal basis of $\LtwoRd$, then $W$ satisfies 
\begin{equation}
\label{Eq:ExpansionW}
W( f ) \stackrel{a.s.}{=} \sum_{n \in \mathbb{N}} W(e_{n}) \left( f , e_{n}\right)_{\LtwoRd}, \quad \forall f \in \LtwoRd,
\end{equation}
where the series is taken in the sense of $\LtwoOmega$ and the random variables $(W(e_{n}))_{n \in \mathbb{N}}$ are pairwise uncorrelated with unit variance. Conversely, if a sequence of uncorrelated with unit variance real random variables $(\epsilon_{n})_{n \in \mathbb{N}}$ is provided,  and if $(e_{n})_{n \in \mathbb{N}}$ is an arbitrary but fixed (real) orthonormal basis of $\LtwoRd$, then the process $(W(f))_{f \in \LtwoRd}$ defined through
\begin{equation}
\label{Eq:DefWfromBasis}
W( f ) = \sum_{n \in \mathbb{N}} \epsilon_{n} \left( f , e_{n}\right)_{\LtwoRd}, \quad \forall f \in \LtwoRd,
\end{equation}
is a real White Noise over $\Rd$. If $\mathcal{L} : \SchwartzRd \to \LtwoRd$ is linear and continuous, then its adjoint $\mathcal{L}^{*} : \LtwoRd \ (=\LtwoRd' \ ) \to \TemperedRd$ can be applied to a White Noise following Proposition \ref{Prop:ExtenstionZtoU} through
\begin{equation}
\label{Eq:L*W}
\mathcal{L}^{*}W (f) := W( \mathcal{L}f ), \quad f \in \SchwartzRd.
\end{equation}
Of course, every linear operator over tempered distributions defined through an adjoint is applicable to a White Noise. 

\section{Linear transformation \textit{from} a White Noise}
\label{Sec:RepresentationZ=LW}

In this section we prove that any GeSP can be expressed as a linear transformation of a White Noise. While the following result requires some technical details, it can be conceived in a quite intuitive manner. Given a GeSP $Z$, we construct a White Noise from which we can obtain $Z$ linearly. We start by applying an enough-high order negative Laplacian operator to $Z$ in order to obtain a mean-square continuous stochastic process. This regular process is then identified to a process indexed by a suitable Hilbert space, for which we can apply the Karhunen-Loève Theorem. The use of the orthonormal basis of this Hilbert space and of the sequence of uncorrelated random variables obtained allows to construct a White Noise from which we can linearly re-obtain $Z$.

\begin{Theo}
\label{Theo:Z=LW}
Let $Z$ be a real GeSP over $\Rd$. Then, there exists a linear and weak-continuous operator $\mathcal{L} : L^{2}(\Rd)\ (=L^{2}(\Rd)'\ ) \to \TemperedRd$ and a real White Noise $W$ over $\Rd$ such that
\begin{equation}
\label{Eq:Z=LW}
Z = \mathcal{L}W.
\end{equation} 
\end{Theo}
\textbf{Proof:} Following regularity Theorem \ref{Prop:RegularityGeSP}, there exists $N \in \mathbb{N}$ such that 
\begin{equation}
\label{Eq:DefY}
Y = (1 - \Delta)^{-\frac{N}{2}}Z
\end{equation}
is (identified with a) mean-square continuous stochastic process over $\Rd$ with covariance function $C_{Y}$ in $C_{PB}(\RdxRd)$. There exists then $M \in \mathbb{N}$, and $C > 0 $ such that
\begin{equation}
\label{Eq:CYBounded}
|C_{Y}(x,y)|\leq C(1+|x^{2}|)^{\frac{M}{2}}(1+|y^{2}|)^{\frac{M}{2}}, \quad \forall (x,y) \in \RdxRd.
\end{equation}
Let us define the positive measure with density over $\Rd$:
\begin{equation}
\label{Eq:DefMu}
d\mu(x) = \frac{dx}{(1+|x|^{2})^{M + \frac{d+1}{2}}},
\end{equation}
which is a finite measure for every $M \in \mathbb{N}$. We consider the complex Hilbert space $\LtwoRdmu$ which is separable and such that $\SchwartzRd \subset L^{2}(\Rd,\mu)$ densely and continuously embedded (the density and continuous embedding of $\SchwartzRd$ in $L^{2}(\Rd)$ can be used to prove this). From the boundedness condition \eqref{Eq:CYBounded}  and the finiteness of $\mu$, it follows that the function $C_{Y}$ is in $L^{2}(\RdxRd , \mu \otimes \mu)$. If $\varphi,\phi \in \SchwartzRd$, then, since $\varphi \otimes \phi \in L^{2}(\RdxRd , \mu \otimes \mu)$, one has
\begin{equation}
\label{Eq:CYcontinuousL2muxL2mu}
\begin{aligned}
|\langle C_{Y}^{*} , \varphi \otimes \phi \rangle | &\leq  \| C_{Y} \|_{L^2(\RdxRd, \mu\otimes\mu)} \| \varphi \otimes \phi \|_{L^2(\RdxRd, \mu\otimes\mu)} \\
 &=  \| C_{Y} \|_{L^2(\RdxRd, \mu\otimes\mu)}  \| \varphi \|_{\LtwoRdmu}\| \phi \|_{\LtwoRdmu}.
\end{aligned}
\end{equation}
From this inequality it can be concluded that the distribution $C_{Y}^{*}$ is in $( L^{2}(\Rd,\mu)\otimes_{\pi}L^{2}(\Rd,\mu) )'$. Hence, following Proposition \ref{Prop:ExtenstionZtoU}, $Y$ can be extended continuously to $\LtwoRdmu$. The stochastic process $(Y(f))_{f \in \LtwoRdmu}$ satisfies hence conditions \ref{Itm:1KL} and \ref{Itm:2KL} of Theorem \ref{Theo:KL}. The covariance Kernel of $Y$ if given by
\begin{equation}
\label{Def:CYasIntegral}
\Cov(Y(f),Y(g)) = \int_{\RdxRd} C_{Y}(x,y)f(x)\overline{g(y)} d(\mu\otimes\mu)(x,y), \quad \forall f,g \in \LtwoRdmu,
\end{equation}
from which we can identify the covariance operator of $Y$, $Q_{Y} : \LtwoRdmu \to \LtwoRdmu $, given by
\begin{equation}
\label{Eq:QY}
Q_{Y}(f) = \int_{\Rd}C_{Y}( \cdot , y)f(y)d\mu(y). 
\end{equation}
$Q_{Y}$ is an integral positive semidefinite operator with kernel $C_{Y} \in L^{2}(\Rd\times\Rd, \mu\otimes\mu )$, which is in addition continuous and for which the function $x \mapsto C_{Y}(x,x) $ is in $L^{1}(\Rd , \mu)$. We conclude that $Q_{Y}$ is a trace-class operator \citep[Theorem 4.3]{brislawn1991traceable}, hence it is in particular compact. Therefore, condition \ref{Itm:3KL} in Theorem \ref{Theo:KL} is also satisfied. We can construct hence the Karhunen-Loève expansion of $Y$ with respect to $L^{2}(\Rd,\mu)$:
\begin{equation}
\label{Eq:KLY}
Y = \sum_{n \in \mathbb{N}} Y(f_{n})\overline{f_{n}}^{*},
\end{equation}
where $(f_{n})_{n \in \mathbb{N}}$ is the orthonormal basis of $L^{2}(\Rd,\mu)$ consisting of the eigenfunctions of $Q_{Y}$. The random variables $(Y(f_{n}))_{n \in \mathbb{N}}$ are pairwise uncorrelated  with variances $\Var(Y(f_{n})) = \lambda_{n} $, being $\lambda_{n}$ the eigenvalue of $Q_{Y}$ associated to the eigenfunction $f_{n}$. The sequence $(\lambda_{n})_{n \in \mathbb{N}}$ is non-negative and in $\ell^{1}(\mathbb{N})$. Let us also define the set
\begin{equation}
\label{Eq:DefN0}
N_{0} := \lbrace n \in \mathbb{N} \ \big| \ \lambda_{n} = 0 \rbrace.
\end{equation}
Let us now consider the operator $\mathcal{L}_{M} : L^{2}(\Rd) \to L^{2}(\Rd, \mu) $ defined by the multiplication by the function $(1+|x|^{2})^{\frac{M}{2}+ \frac{d+1}{4}}$:
\begin{equation}
\label{Eq:DefLM}
\mathcal{L}_{M}( f ) = (1+|x|^{2})^{\frac{M}{2}+ \frac{d+1}{4}} f.
\end{equation}
It is easy to verify that this operator is linear, continuous and bijective, and its inverse, $\mathcal{L}_{M}^{-1}$, is the multiplication by the function $(1+|x|^{2})^{-\frac{M}{2}- \frac{d+1}{4}}$. In addition, $\mathcal{L}_{M}$ is the multiplication with a multiplicator of the Schwartz space (a smooth function with polynomially bounded derivatives, see \citet[Example 7, Section V.3]{reed1980methods} or \cite[Definition 25.3]{treves1967topological}), hence it is a continuous linear operator $\SchwartzRd \to \SchwartzRd$, which is in addition bijective (its inverse is also a multiplication by a multiplicator of the Schwartz space). The adjoint $\mathcal{L}_{M}^{*}$ \textit{is also} the multiplication by the function $(1+|x|^{2})^{\frac{M}{2}}$ applied over $\TemperedRd$, and its a bijective endomorphism. The same goes for $(\mathcal{L}^{*}_{M})^{-1}$. Hence, both operators can be applied to any GeSP without difficulties. We define then the functions
\begin{equation}
\label{Eq:Defgn}
g_{n} =  \mathcal{L}_{M}^{-1}(f_{n}) = (1+|x|^{2})^{-\frac{M}{2} - \frac{d+1}{4}}f_{n}, \quad  \forall n \in \mathbb{N}.
\end{equation}
Since $(f_{n})_{n \in \mathbb{N}}$ is an orthonormal basis of $\LtwoRdmu$, it can be easily concluded that $(g_{n})_{n \in \mathbb{N}}$ is an orthonormal basis of $\LtwoRd$. Now, we consider a collection of random variables with unit variance $(\epsilon_{n})_{n \in N_{0}}$, pairwise uncorrelated and independent of $Z$ (and hence of $Y$)\footnote{At this point, we suppose that the probability space $\OmegaAP$ is big enough so these random variables can be constructed. If it is not, since we are just adding independent random variables, it is enough to consider a new probability space by taking a suitable product of probability spaces, as usual.}. We construct then the GeSP over $\Rd$:
\begin{equation}
\label{Eq:DefWsuchthatZ=LW}
W := \sum_{n \in \mathbb{N}\setminus N_{0}} \frac{Y(f_{n})}{\sqrt{\lambda_{n}}}g_{n}^{*}  + \sum_{n \in N_{0}} \epsilon_{n}g_{n}^{*}.
\end{equation}
Since the random variables $\frac{Y(f_{n})}{\sqrt{\lambda_{n}}}$ and $\epsilon_{n}$ are pairwise uncorrelated with unit variance, the process $W$ is a White Noise over $\Rd$ (Eq. \eqref{Eq:DefWfromBasis}).

Now we consider the operator, $\mathcal{L}_{Q}^{*} : \LtwoRd' \to \TemperedRd$ given by\footnote{While the functions $g_{n}$ are real, we decided to make explicit this construction considering the complex conjugates $\overline{g_{n}}$ just in case that some extension to the case of complex GeSP is desired.}
\begin{equation}
\label{Eq:DefLQ*}
\mathcal{L}_{Q}^{*}(\cdot) = \sum_{n \in \mathbb{N}} \sqrt{\lambda_{n}} \langle \cdot ,  \overline{g_{n}} \rangle g_{n}^{*}.
\end{equation}
It is left to the reader to verify that $\mathcal{L}_{Q}^{*}$ is the adjoint of the operator $\mathcal{L}_{Q} : \SchwartzRd \to \LtwoRd$ given by
\begin{equation}
\mathcal{L}_{Q}(\varphi) = \sum_{n \in \mathbb{N}} \sqrt{\lambda} ( \varphi , \overline{g_{n}} )_{L^{2}(\Rd)} \overline{g_{n}}, \quad \forall \varphi \in \SchwartzRd.
\end{equation}
Since $(\lambda_{n})_{n \in \mathbb{N}}$ is bounded, $\mathcal{L}_{Q}$ is continuous ($\SchwartzRd$ is continuously embedded in $\LtwoRd$):
\begin{equation}
\label{Eq:LQcontinuous}
\| \mathcal{L}_{Q}(\varphi) \|_{\LtwoRd}^{2} = \sum_{n \in \mathbb{N}} |(\varphi,\overline{g_{n}})_{\LtwoRd}|^{2}\lambda_{n} \leq \max_{n \in \mathbb{N}}\lambda_{n} \| \varphi \|_{\LtwoRd}^{2}, \quad \forall \varphi \in \SchwartzRd.
\end{equation}
Hence, $\mathcal{L}_{Q}^{*}$ can be applied to $W$. We also have $\mathcal{L}_{Q}^{*}(g_{n}^{*}) = \sqrt{\lambda_{n}}g_{n}^{*}$ for every $n \in \mathbb{N}$.

 We finally consider the operator $\mathcal{L} : \LtwoRd \to \TemperedRd$ given by
\begin{equation}
\label{Eq:DefL}
\mathcal{L} = (1 - \Delta)^{\frac{N}{2}} \circ \mathcal{L}_{M}^{*} \circ \mathcal{L}_{Q}^{*},
\end{equation}
which is the adjoint of a continuous linear operator from $L^{2}(\Rd)$ to $\TemperedRd$. We obtain thus, interpreting the infinite sums in a mean-square-weak-$*$ sense, 
\begin{equation}
\label{Eq:ProofZ=LW}
\begin{aligned}
\mathcal{L}W &=  (1-\Delta)^{\frac{N}{2}} \mathcal{L}_{M}^{*}\mathcal{L}_{Q}^{*} \left( \sum_{n \in \mathbb{N}\setminus N_{0}} \frac{Y(f_{n})}{\sqrt{\lambda_{n}}} g_{n}^{*} + \sum_{n \in N_{0}} \epsilon_{n}g_{n}^{*} \right) \\
&= (1-\Delta)^{\frac{N}{2}}\mathcal{L}_{M}^{*}\left( \sum_{n \in \mathbb{N}\setminus N_{0}}Y(f_{n})g_{n}^{*}  \right)  &( \mathcal{L}_{Q}^{*} ( g_{n}^{*} ) = 0 \hbox{ if } n \in N_{0} ) \\
&= (1-\Delta)^{\frac{N}{2}}\left( \sum_{n \in \mathbb{N}\setminus N_{0}} Y(f_{n})f_{n}^{*} \right)  &( \mathcal{L}_{M}^{*}(g_{n}^{*} ) = f_{n}^{*} ) \\
&= (1-\Delta)^{\frac{N}{2}}\left( \sum_{n \in \mathbb{N} } Y(f_{n})f_{n}^{*} \right)  &(Y(f_{n}) \stackrel{a.s}{=} 0 \hbox{ if } n \in N_{0} ) \\
&=(1-\Delta)^{\frac{N}{2}}Y = Z. \quad \blacksquare
\end{aligned}
\end{equation}

In the following, we will freely use the objects defined in the proof of Theorem \ref{Theo:Z=LW}.

One important Corollary obtained from the proof of Theorem \ref{Theo:Z=LW} is the orthogonal decomposition of a GeSP. This is, in some sense, a generalization of the Karhunen-Loève expansion to the case of random distributions. It also implies that any GeSP can be constructed from an at most countable quantity of uncorrelated random variables with non-null variance. It also allows to obtain a similar representation for the GeSP once an applicable linear operator is applied.

\begin{Corol}[\textbf{Orthogonal representation of a GeSP}]
\label{Corol:Zexpression}
Let $Z$ be a real GeSP over $\Rd$. Then, there exists a sequence of pairwise uncorrelated real random variables $(Z_{n})_{n \in \mathbb{N}}$ and a sequence of linearly independent real tempered distributions $(T_{n})_{n \in \mathbb{N}}$ such that
\begin{equation}
\label{Eq:ZexpandedTn}
Z = \sum_{n \in \mathbb{N}}Z_{n} T_{n},
\end{equation}
where the series is interpreted in a mean-square-$\TemperedRd$-weak sense. The covariance distribution of $Z$ is then given by
\begin{equation}
\label{Eq:CZexpanded}
C_{Z} = \sum_{n \in \mathbb{N}} \lambda_{n} T_{n} \otimes T_{n},
\end{equation}
with $\lambda_{n}$ non-negative numbers such that $\sum_{n \in \mathbb{N}}\lambda_{n} < \infty $. The series \eqref{Eq:CZexpanded} is interpreted in a $(\SchwartzRd\otimes_{\pi}\SchwartzRd)'$-weak sense. In addition, if $\mathscr{U}(\Rd)$ is a space as in Proposition \ref{Prop:ExtenstionZtoU}, then for every $n \in \mathbb{N}\setminus N_{0}$ the distribution $T_{n}$ is in $\mathscr{U}'(\Rd)$, the convergence of the series \eqref{Eq:ZexpandedTn} is also in the mean-square-$\mathscr{U}'(\Rd)$-weak sense, and that one of the series \eqref{Eq:CZexpanded} is in the $(\mathscr{U}(\Rd)\otimes_{\pi}\mathscr{U}(\Rd))'$-weak sense. Hence, if $\mathcal{L} : \SchwartzRd \to \mathscr{U}(\Rd)$ is linear and continuous, its adjoint $\mathcal{L}^{*}$ is applicable to $Z$ and to every  $T_{n}$, $n \in \mathbb{N}\setminus N_{0}$. We have in such case,
\begin{equation}
\label{Eq:L*Zexpanded}
\mathcal{L}^{*}Z = \sum_{n \in \mathbb{N}\setminus N_{0}} Z_{n} \mathcal{L}^{*}T_{n}
\end{equation}
in a mean-square-$\TemperedRd$-weak sense, and 
\begin{equation}
\label{Eq:CL*Zexpanded}
C_{\mathcal{L}^{*}Z} = \sum_{n\in \mathbb{N}\setminus N_{0}} \lambda_{n} \mathcal{L}^{*}T_{n} \otimes \overline{\mathcal{L}^{*}}T_{n}
\end{equation}
in a $(\SchwartzRd\otimes_{\pi}\SchwartzRd)'$-weak sense. 
\end{Corol}

\textbf{Proof:} Take $Z_{n} = Y(f_{n})$ and $T_{n} = (1-\Delta)^{\frac{N}{2}}\mathcal{L}^{*}_{M}g_{n}^{*}$. The values $\lambda_{n} = \Var(Z_{n})$ are the eigenvalues of $Q_{Y}$ and they are in $\ell^{1}(\mathbb{N})$ since $Q_{Y}$ is trace-class. The mode of convergence of the series \eqref{Eq:CZexpanded} is evident. 

Let $\mathscr{U}(\Rd)$ as in Proposition \ref{Prop:ExtenstionZtoU}. Since $C_{Z} \in (\mathscr{U}(\Rd)\otimes_{\pi}\mathscr{U}(\Rd) )'$, there exists a constant $C > 0$ and a seminorm $p$ over $\mathscr{U}(\Rd)$ such that
\begin{equation}
 \langle C_{Z} , \varphi \otimes \overline{\varphi} \rangle  \leq C p(\varphi)^{2}, \quad \forall \varphi \in \SchwartzRd.
\end{equation}
Hence, if $n \in \mathbb{N}\setminus N_{0}$, one has
\begin{equation}
\label{Eq:TnContinuousU}
\lambda_{n}\left| \langle T_{n} , \varphi \rangle \right|^{2} \leq \sum_{k \in \mathbb{N}\setminus N_{0}} \lambda_{k} | \langle T_{k} , \varphi \rangle |^{2} = \langle C_{Z} , \varphi \otimes \overline{\varphi} \rangle \leq C p(\varphi)^{2}, 
\end{equation}
which implies
\begin{equation}
|\langle T_{n} , \varphi \rangle | \leq \sqrt{\frac{C}{\lambda_{n}}} p(\varphi), \quad \forall \varphi \in \SchwartzRd,
\end{equation}
hence $T_{n} \in \mathscr{U}'(\Rd)$ for $n \in \mathbb{N}\setminus N_{0}$ (we recall that $\SchwartzRd$ is dense in $\mathscr{U}(\Rd)$). Inequalities \eqref{Eq:TnContinuousU} also imply that the series $\sum_{k \in \mathbb{N}\setminus N_{0}} \lambda_{k} | \langle T_{k} , \varphi \rangle |^{2}$ can be extended continuously to $\varphi \in \mathscr{U}(\Rd)$, which implies the new modes of convergence of series \eqref{Eq:ZexpandedTn} and \eqref{Eq:CZexpanded}. We finally obtain the expressions of $\mathcal{L}^{*}Z$ and its covariance through
\begin{equation}
\mathcal{L}^{*}Z(\varphi) = Z(\mathcal{L}\varphi ) = \sum_{n \in \mathbb{N}\setminus N_{0}} Z_{n}\langle T_{n} , \mathcal{L}\varphi \rangle = \sum_{n \in \mathbb{N}\setminus N_{0}} Z_{n}\langle \mathcal{L}^{*}T_{n} , \varphi \rangle, \quad \varphi \in \SchwartzRd.
\end{equation}
Equation \eqref{Eq:CL*Zexpanded} is then obtained directly. $\blacksquare$

We remark that the so-constructed expansion \eqref{Eq:ZexpandedTn} is not unique. Indeed, one could do this kind of construction using any arbitrary $N$ satisfying that $Y$ is mean-square continuous, any arbitrary $M$ such that \eqref{Eq:CYBounded} holds, or with another convenient measure $\mu$.

\section{Linear transformation \textit{into} a White Noise}
\label{Sec:LZ=W}

In this section we look for necessary and sufficient conditions for the existence of an operator which transforms a GeSP into a White Noise. Hence, we basically show the existence of an operator $\mathcal{L}$ and of a White Noise $W$ which, under suitable conditions, makes $Z$ to satisfy the SPDE over $\Rd$ of the form
\begin{equation}
\label{Eq:SPDELZ=W}
\mathcal{L}Z = W.
\end{equation} 

We recall that when we say that an operator $\mathcal{L}$ is \textit{applicable} to $Z$, it is meant in the sense of Proposition \ref{Prop:ExtenstionZtoU}. That is, we assume the existence of a topological vector space $\mathscr{U}(\Rd)$ satisfying the conditions in this Proposition, and that $\mathcal{L}$ is the adjoint of a continuous linear operator from $\SchwartzRd$ to $\mathscr{U}(\Rd)$. In order to construct such a space, we will take advantage from the orthogonal representation from Corollary \ref{Corol:Zexpression} and construct a separable Hilbert space related to this representation. This space is deeply related with the concept of Reproducing Kernel Hilbert space, see \citet{meidan1979reproducing}. While we could use some non-direct abstract constructions (such as taking completitions or using quotient spaces), in this work we shall use a ``\textit{direct and concrete}'' construction of such a space, that is, we will make precise what are its elements.

We begin by defining the bounded sequence of strictly positive numbers
\begin{equation}
\label{Eq:Deflambdatilde}
\tilde{\lambda_{n}} := \begin{cases}
\lambda_{n} & \hbox{ if } n \in \mathbb{N}\setminus N_{0}, \\
1 & \hbox{ if } n \in N_{0}.
\end{cases}
\end{equation}
We consider then the tempered distributions in Corollary \ref{Corol:Zexpression}, $T_{n} = (1-\Delta)^{\frac{N}{2}}\mathcal{L}_{M}^{*}g_{n}^{*}$ for $n \in \mathbb{N}$. We define the following vector space, which is a subspace of the algebraic dual of a subspace of $\TemperedRd$:
\begin{equation}
\label{Eq:DefHZ}
\mathscr{H}_{Z}(\Rd) = \lbrace u \in \Span \lbrace (T_{n})_{n \in \mathbb{N}} \rbrace^{*} \ \big| \  \sum_{n \in \mathbb{N}} \tilde{\lambda_{n}} | \langle u , T_{n} \rangle |^{2} < \infty \rbrace.
\end{equation} 
This space is endowed with the inner product
\begin{equation}
\label{Eq:InnerProductHZ}
( u , v )_{\mathscr{H}_{Z}(\Rd)} =  \sum_{n \in \mathbb{N}} \tilde{\lambda_{n}} \langle u , T_{n} \rangle  \overline{ \langle v , T_{n} \rangle }.
\end{equation}
Using that $\tilde{\lambda_{n}} > 0 $ for every $n$ and the absolute convergence of the involved series, it can be easily concluded that $( \cdot , \cdot )_{\mathscr{H}_{Z}(\Rd)}$ is indeed an inner product over $\mathscr{H}_{Z}(\Rd)$. By considering the isometric and bijective mapping $\mathscr{H}_{Z}(\Rd)\to \ell^{2}(\mathbb{N})$ given by $ u \mapsto ( \sqrt{\tilde{\lambda_{n}}} \langle u , T_{n} \rangle )_{n \in \mathbb{N}}$, it is concluded that $\mathscr{H}_{Z}(\Rd)$ is a separable Hilbert space.

Let us identify $\SchwartzRd$ with a subspace of $\mathscr{H}_{Z}(\Rd)$. Since $T_{n} \in \TemperedRd$, every $\varphi \in \SchwartzRd$ will be canonically identified with the element $\varphi^{**} \in \Span \lbrace (T_{n})_{n \in \mathbb{N}} \rbrace^{*} $ given by
\begin{equation}
\label{Eq:DefVarphi**}
\langle \varphi^{**} , T_{n} \rangle := \langle T_{n} , \varphi \rangle, \quad \forall n \in \mathbb{N}.
\end{equation}
Let us verify that $\varphi^{**}$ is indeed in $\mathscr{H}_{Z}(\Rd)$. For that, we argue that 

\begin{equation}
\label{Eq:varphi**inHZ}
\begin{aligned}
\| \varphi^{**} \|_{\mathscr{H}_{Z}(\Rd)}^{2} &= \sum_{n \in \mathbb{N}}\tilde{\lambda_{n}}| \langle \varphi^{**} , T_{n} \rangle |^{2} \\
&= \sum_{n \in \mathbb{N}}\tilde{\lambda_{n}}| \langle  T_{n} , \varphi \rangle |^{2} \\
&= \sum_{n \in \mathbb{N}}\tilde{\lambda_{n}}| \langle (1-\Delta)^{\frac{N}{2}}\mathcal{L}^{*}_{M}g_{n}^{*} , \varphi \rangle |^{2}  \\
&= \sum_{n \in \mathbb{N}}\tilde{\lambda_{n}} \left| \left( \mathcal{L}_{M}(1-\Delta)^{\frac{N}{2}}\varphi , \overline{g_{n}} \right)_{L^{2}(\Rd)} \right|^{2} \\
&\leq \sup_{n \in \mathbb{N}} \tilde{\lambda_{n}} \| \mathcal{L}_{M}(1-\Delta)^{\frac{N}{2}}\varphi \|^{2}_{L^{2}(\Rd)} < \infty,
\end{aligned}
\end{equation}
where we have used that $(g_{n})_{n \in \mathbb{N}}$ is an orthonormal basis of $L^{2}(\Rd)$ and that $\mathcal{L}_{M}(1-\Delta)^{\frac{N}{2}}\varphi \in \SchwartzRd \subset L^{2}(\Rd)$. Hence $\varphi^{**} \in \mathscr{H}_{Z}(\Rd)$.

$\SchwartzRd$ ``is'' not the only interesting subspace of $\mathscr{H}_{Z}(\Rd)$. Let us consider the following subspace of tempered distributions, which is defined similarly to Sobolev spaces:
\begin{equation}
(1-\Delta)^{-\frac{N}{2}}\mathcal{L}_{M}^{-1} L^{2}(\Rd) = \lbrace T \in \TemperedRd \ \big| \ \exists f \in L^{2}(\Rd) \hbox{ s.t. } T = (1-\Delta)^{-\frac{N}{2}}\mathcal{L}_{M}^{-1}f \rbrace.
\end{equation}
Let $T \in (1-\Delta)^{-\frac{N}{2}}\mathcal{L}_{M}^{-1} L^{2}(\Rd)$ and let $f \in L^{2}(\Rd)$ such that $T = (1-\Delta)^{-\frac{N}{2}}\mathcal{L}_{M}^{-1}$. Then, we define its canonically associated element $T^{**} \in \mathscr{H}_{Z}(\Rd)$ as
\begin{equation}
\label{Eq:DefT**}
\langle T^{**} , T_{n} \rangle := \left( f ,  \overline{g_{n}} \right)_{L^{2}(\Rd)}, \quad \forall n \in \mathbb{N}.
\end{equation}
It is clear that $\| T^{**} \|_{\mathscr{H}_{Z}(\Rd)}^{2} \leq \sup_{n \in \mathbb{N}}\tilde{\lambda_{n}} \| f \|^{2}_{L^{2}(\Rd)} < \infty$, hence $T^{**}$ is indeed in $\mathscr{H}_{Z}(\Rd)$. This construction may seem artificial, but it actually follows the intuitive rule of the application of the inverse transpose operator:
\begin{equation}
\langle T^{**} , T_{n} \rangle = ``\langle (1-\Delta)^{-\frac{N}{2}}\mathcal{L}_{M}^{-1}f , (1-\Delta)^{\frac{N}{2}}\mathcal{L}_{M}^{*}g_{n}^{*} \rangle" = \langle f^{*} , g_{n} \rangle = (f , \overline{g_{n}} )_{L^{2}(\Rd)}.
\end{equation}
In addition, if $\varphi \in \SchwartzRd$ since its associated distribution $\varphi^{*} \in \TemperedRd$ clearly is in $(1-\Delta)^{-\frac{N}{2}}\mathcal{L}_{M}^{-1} L^{2}(\Rd)$,  we have
\begin{equation}
\begin{aligned}
\langle [\varphi^{*}]^{**} , T_{n} \rangle &= \left( \mathcal{L}_{M}(1-\Delta)^{\frac{N}{2}}\varphi , \overline{g_{n}} \right)_{L^{2}(\Rd)} \\
&= \langle  g_{n}^{*} , \mathcal{L}_{M}(1-\Delta)^{\frac{N}{2}}\varphi  \rangle \\
&= \langle  (1-\Delta)^{\frac{N}{2}}\mathcal{L}^{*}_{M}g_{n}^{*} , \varphi  \rangle \\
&= \langle T_{n} , \varphi \rangle \\
&= \langle \varphi^{**} , T_{n} \rangle, \quad \forall n \in \mathbb{N}.
\end{aligned}
\end{equation}
Hence $[\varphi^{*}]^{**} = \varphi^{**}$, which justifies the canonical definition of $T^{**}$.

Let us now consider the distributions in $(1-\Delta)^{-\frac{N}{2}}\mathcal{L}_{M}^{-1} L^{2}(\Rd)$ defined through
\begin{equation}
S_{n} := \frac{(1-\Delta)^{-\frac{N}{2}}\mathcal{L}_{M}^{-1}\overline{g_{n}}}{\sqrt{\tilde{\lambda_{n}}}}, \quad n \in \mathbb{N}.
\end{equation}
Then,
\begin{equation}
\begin{aligned}
\left( S_{j}^{**} , S_{k}^{**}\right)_{\mathscr{H}_{Z}(\Rd)} &= \sum_{n \in \mathbb{N}} \tilde{\lambda_{n}} \langle S_{j}^{**} , T_{n} \rangle \overline{\langle S_{k}^{**} , T_{n} \rangle} \\
&= \sum_{n \in \mathbb{N}} \tilde{\lambda_{n}} \left( \frac{\overline{g_{j}}}{\sqrt{\tilde{\lambda_{j}}}}  , \overline{g_{n}}\right)_{L^{2}(\Rd)}\overline{ \left(   \frac{\overline{g_{k}}}{\sqrt{\tilde{\lambda_{k}}}}  , \overline{g_{n}} \right)_{L^{2}(\Rd)}} \\
&= \sum_{n \in \mathbb{N}} \frac{\tilde{\lambda_{n}}}{\sqrt{\tilde{\lambda_{j}}\tilde{\lambda_{k}}}} \delta_{j,n} \delta_{n,k} = \delta_{j,k}, \quad \forall j,k  \in \mathbb{N}.
\end{aligned}
\end{equation}
We conclude that the elements $(S_{n}^{**})_{n \in \mathbb{N}}$ form an orthonormal system of $\mathscr{H}_{Z}(\Rd)$. It is in addition a complete basis, which can be verified using the separating points criterion:
\begin{equation}
\begin{aligned}
\left( u , S_{j}^{**} \right)_{\mathscr{H}_{Z}(\Rd)} = 0 \quad \forall j \in \mathbb{N} &\Leftrightarrow \sum_{n \in \mathbb{N}} \tilde{\lambda_{n}} \langle u , T_{n} \rangle \overline{\langle S_{j}^{**} , T_{n} \rangle} =  0 \quad \forall j \in \mathbb{N} \\
&\Leftrightarrow \sum_{n \in \mathbb{N}} \sqrt{\tilde{\lambda_{n}}} \langle u , T_{n} \rangle \delta_{j,n} = 0 \quad \forall j \in \mathbb{N}, \\
&\Leftrightarrow \sqrt{\tilde{\lambda_{j}}} \langle u , T_{j} \rangle = 0 \quad \forall j \in \mathbb{N}, \\
& \Leftrightarrow u = 0 \quad \hbox{(since $\tilde{\lambda_{j}} > 0 $)}.
\end{aligned}
\end{equation}

Let us now consider the dual space $\mathscr{H}_{Z}'(\Rd)$. By Riesz Representation, every  member $T \in \mathscr{H}_{Z}'(\Rd)$ is canonically and isometrically identified with a unique member $u_{T} \in \mathscr{H}_{Z}(\Rd)$ following
\begin{equation}
\langle T , u \rangle = ( u , u_{T} )_{\mathscr{H}_{Z}(\Rd)}, \quad \forall u \in \mathscr{H}_{Z}(\Rd).
\end{equation}
Of course, this isometry is considered in the sense of the inner product in $\mathscr{H}_{Z}'(\Rd)$ given by
\begin{equation}
\left( T , S \right)_{\mathscr{H}_{Z}'(\Rd)} := \left( u_{T} , u_{S} \right)_{\mathscr{H}_{Z}(\Rd)}, \quad \forall T,S \in \mathscr{H}_{Z}'(\Rd).
\end{equation}
As known, the mapping $T \mapsto u_{T}$ maps orthonormal bases in orthonormal bases. Hence, we can consider the elements $(h_{n})_{n \in \mathbb{N}} \in \mathscr{H}_{Z}'(\Rd)$ defined through
\begin{equation}
\langle h_{n} , u \rangle := \left( u , S_{n}^{**} \right)_{\mathscr{H}_{Z}(\Rd)}, \quad \forall u \in \mathscr{H}_{Z}(\Rd), \forall n \in \mathbb{N},
\end{equation}
and conclude that they form a complete orthonormal basis of $\mathscr{H}_{Z}'(\Rd)$. Let us see to which tempered distribution these elements could be canonically identified. That is, we need to verify for every $j \in \mathbb{N}$ if there exists $T \in \TemperedRd$ such that
\begin{equation}
\langle h_{j} , \varphi^{**} \rangle = \langle T , \varphi \rangle, \quad \forall \varphi \in \SchwartzRd.
\end{equation}
To do this, we remark that
\begin{equation}
\label{Eq:hn=sqrt(ln)Tnargument}
\begin{aligned}
\langle h_{j} , \varphi^{**} \rangle &= \left( \varphi^{**} , S_{j}^{**} \right)_{\mathscr{H}_{Z}(\Rd)} \\
&= \sum_{n \in \mathbb{N}} \tilde{\lambda_{n}}\langle \varphi^{**} , T_{n} \rangle \overline{ \langle  S_{j}^{**} , T_{n} \rangle} \\
&= \sum_{n \in \mathbb{N}}\tilde{\lambda_{n}} \langle T_{j} , \varphi \rangle\frac{\delta_{j,n}}{\sqrt{\tilde{\lambda_{j}}}} \\
&= \sqrt{\tilde{\lambda_{j}}}\langle T_{j} , \varphi \rangle = \langle \sqrt{\tilde{\lambda_{j}}}T_{j} , \varphi \rangle , \quad \forall \varphi \in \SchwartzRd.
\end{aligned} 
\end{equation}
We conclude that 
\begin{equation}
\label{Eq:hn=sqrt(ln)Tn}
h_{n} \simeq \sqrt{\tilde{\lambda_{n}}}T_{n}, \quad \forall n \in \mathbb{N},
\end{equation}
where $\simeq$ means that $h_{n}$ is canonically identified with $\sqrt{\tilde{\lambda_{n}}}T_{n}$ in the sense of equation \eqref{Eq:hn=sqrt(ln)Tnargument}. The distributions $(\sqrt{\tilde{\lambda_{n}}}T_{n})_{n \in \mathbb{N}}$ ``form'' then an orthonormal basis of $\mathscr{H}_{Z}'(\Rd)$.

\begin{Prop}
$\SchwartzRd \subset \mathscr{H}_{Z}(\Rd)$ densily and continuously embedded, and in addition $C_{Z} \in (\mathscr{H}_{Z}(\Rd)\otimes_{\pi}\mathscr{H}_{Z}(\Rd))'$.
\end{Prop}

\textbf{Proof:} The continuous immersion of $\SchwartzRd$ in $\mathscr{H}_{Z}(\Rd)$ comes simply from the inequality \eqref{Eq:varphi**inHZ} and from the continuity of $\mathcal{L}_{M}(1-\Delta)^{\frac{N}{2}}$ as an operator from $\SchwartzRd $ to $\SchwartzRd$. Indeed, since this operator is continuous, and $\SchwartzRd$ is continuously embedded in $L^{2}(\Rd)$, there exists $C' > 0 $ and a seminorm of the Schwartz space $\mathcal{N}_{m}$ such that
\begin{equation}
\| \varphi^{**} \|_{\mathscr{H}_{Z}(\Rd)} \leq \sqrt{\sup_{n \in \mathbb{N}}\tilde{\lambda_{n}}} \| \mathcal{L}_{M}(1-\Delta)^{\frac{N}{2}}\varphi \|_{L^{2}(\Rd)} \leq  \sqrt{\sup_{n \in \mathbb{N}}\tilde{\lambda_{n}}} C' \mathcal{N}_{m}(\varphi), \quad \forall \varphi \in \SchwartzRd. 
\end{equation}
The density can be argued as follows. Since the distributions $(S_{n})_{n \in \mathbb{N}}$ ``form'' an orthonormal basis of $\mathscr{H}_{Z}(\Rd)$, one has that $\Span\lbrace(S_{n})_{n \in \mathbb{N}} \rbrace \subset (1-\Delta)^{-\frac{N}{2}}\mathcal{L}_{M}^{-1}L^{2}(\Rd)$ ``is'' dense in $\mathscr{H}_{Z}(\Rd)$, hence the space $(1-\Delta)^{-\frac{N}{2}}\mathcal{L}_{M}^{-1}L^{2}(\Rd)$ ``is'' dense in $\mathscr{H}_{Z}(\Rd)$. It suffices then to check that to every $T \in (1-\Delta)^{-\frac{N}{2}}\mathcal{L}_{M}^{-1}L^{2}(\Rd)$ there is $\varphi \in \SchwartzRd$ such that $\varphi^{**}$ is arbitrarily close to $T^{**}$ according to the norm of $\mathscr{H}_{Z}(\Rd)$. Let us then take $f \in L^{2}(\Rd)$ such that $T = (1-\Delta)^{-\frac{N}{2}}\mathcal{L}_{M}^{-1}f$. Since $\SchwartzRd$ is dense in $L^{2}(\Rd)$, there exists $\phi \in \SchwartzRd$ arbitrarily close to $f$ in the norm of $L^{2}(\Rd)$. Let us take $\varphi = (1-\Delta)^{-\frac{N}{2}}\mathcal{L}_{M}^{-1}\phi $. We have then,
\begin{equation}
\| T^{**} - \varphi^{**} \|_{\mathscr{H}_{Z}(\Rd)}^{2} \leq \sup_{n \in \mathbb{N}} \tilde{\lambda_{n}} \| f - \phi \|_{L^{2}(\Rd)}^{2},
\end{equation}
from which the density of $\SchwartzRd$ in $\mathscr{H}_{Z}(\Rd)$ is concluded.

We finally have, using the representation \eqref{Eq:CZexpanded} of the covariance $C_{Z}$,
\begin{equation}
\begin{aligned}
\langle C_{Z} , \varphi \otimes \overline{\varphi} \rangle &= \sum_{n \in \mathbb{N}\setminus N_{0}}\lambda_{n} | \langle T_{n} , \varphi \rangle |^{2} \\
&=\sum_{n \in \mathbb{N}\setminus N_{0}} | \langle \sqrt{\lambda_{n}} T_{n} , \varphi \rangle |^{2} \\
&\leq \sum_{n \in \mathbb{N}} | \langle h_{n} , \varphi^{**} \rangle |^{2} \\
&= \sum_{n \in \mathbb{N}} \left| \left( \varphi^{**} , S_{n}^{**} \right)_{\mathscr{H}_{Z}(\Rd)} \right|^{2} \\
&= \| \varphi^{**} \|^{2}_{\mathscr{H}_{Z}(\Rd)}.
\end{aligned} 
\end{equation}
Hence, by Corollary \ref{Corol:UHilbert}, $C_{Z} \in (\mathscr{H}_{Z}(\Rd)\otimes_{\pi}\mathscr{H}_{Z}(\Rd))'$. $\blacksquare$

We have hence constructed a separable Hilbert space $\mathscr{H}_{Z}(\Rd)$ satisfying conditions in Proposition \ref{Prop:ExtenstionZtoU} and we have studied some of its properties. Now, using this space, we answer to the question about when we can transform $Z$ into a White Noise linearly. This will simply come from the possibility of constructing a linear and continuous bijection between $L^{2}(\Rd)$ and a subspace of $\mathscr{H}_{Z}(\Rd)$, which hence must be infinite dimensional.

\begin{Theo}
\label{Theo:LZ=W}
Let $Z$ be a real GeSP over $\Rd$. Then, there exists a linear operator $\mathcal{L}$ applicable to $Z$ such that $\mathcal{L}Z$ is a White Noise if and only if the set $\mathbb{N}\setminus N_{0}$ is infinite.
\end{Theo}

We remark that the condition of $\mathbb{N}\setminus N_{0}$ being infinite is equivalent to $Q_{Y}$ having infinite rank. This tell us that the series \eqref{Eq:ZexpandedTn} and \eqref{Eq:CZexpanded} are actually infinite series and not finite sums.

\textbf{Proof: } If $\mathbb{N}\setminus N_{0}$ is infinite, then there exists a bijection $\gamma : \mathbb{N}\setminus N_{0} \to \mathbb{N}$. Using that $(g_{n})_{n \in \mathbb{N}}$ is an orthonormal basis of $L^{2}(\Rd)$\footnote{We could have use any orthonormal basis of $L^{2}(\Rd)$. The choice of $(g_{n})_{n \in \mathbb{N}}$ is only done in order to save notation.}, let us then consider the linear and continuous operator $\mathcal{L}_{\gamma} : L^{2}(\Rd) \to \mathscr{H}_{Z}(\Rd)$ constructed such as
\begin{equation}
\mathcal{L}_{\gamma}(g_{n}) = S^{**}_{\gamma^{-1}(n)}, \quad \forall n \in \mathbb{N}.
\end{equation}
Hence, the action of $\mathcal{L}_{\gamma}$ over any $f \in L^{2}(\Rd)$ is expressed through the decomposition of $f$ in the basis $(g_{n})_{n \in \mathbb{N}}$:
\begin{equation}
\mathcal{L}_{\gamma}(f) = \mathcal{L}_{\gamma}\left( \sum_{n \in \mathbb{N}} \left( f , g_{n}\right)_{L^{2}(\Rd)} g_{n} \right) = \sum_{n \in \mathbb{N}} \left( f , g_{n}\right)_{L^{2}(\Rd)} S_{\gamma^{-1}(n)}^{**}.
\end{equation}
Since $(S_{\gamma^{-1}(n)}^{**})_{n \in \mathbb{N}}$ is an orthonormal set of $\mathscr{H}_{Z}(\Rd)$, it is easy to verify that $\| \mathcal{L}_{\gamma}(f) \|_{\mathscr{H}_{Z}(\Rd)} = \| f \|_{L^{2}(\Rd)}$, from which it is clear that $\mathcal{\gamma}$ is indeed continuous. In consequence, its adjoint $\mathcal{L}_{\gamma}^{*} : \mathscr{H}_{Z}'(\Rd) \to L^{2}(\Rd)' = L^{2}(\Rd) \subset \TemperedRd$ is weak-$*$ continuous and it can thus be applied to $Z$. Let us consider its action over the elements $h_{k} \simeq  \sqrt{\tilde{\lambda_{k}}}T_{k} \in \TemperedRd$:
\begin{equation}
\begin{aligned}
\langle \mathcal{L}_{\gamma}^{*}(h_{k}) , \varphi \rangle &= \langle h_{k} , \mathcal{L}_{\gamma}(\varphi) \rangle \\
&= \left( \mathcal{L}_{\gamma}(\varphi) , S_{k}^{**} \right)_{\mathscr{H}_{Z}(\Rd)} \\
&= \left( \sum_{n \in \mathbb{N}} (\varphi , g_{n} )_{L^{2}(\Rd)}S_{\gamma^{-1}(n)}^{**} , S_{k}^{**}  \right)_{\mathscr{H}_{Z}(\Rd)} \\
&= \sum_{n \in \mathbb{N}} (\varphi , g_{n} )_{L^{2}(\Rd)} \delta_{\gamma^{-1}(n) , k} \\
&= ( \varphi , g_{\gamma(n)} )_{L^{2}(\Rd)} = \langle \overline{g_{\gamma(n)}^{*}} , \varphi \rangle, \quad \forall \varphi \in \SchwartzRd.
\end{aligned}
\end{equation}
Consequently, since the functions $(g_{n})_{n \in \mathbb{N}}$ are real,
\begin{equation}
\mathcal{\gamma}^{*}( \underbrace{h_{k}}_{\simeq \sqrt{\tilde{\lambda_{k}}}T_{k} } ) = g_{\gamma(n)}^{*}, \quad \forall k \in \mathbb{N}.
\end{equation}
Following then the expansion \eqref{Eq:CL*Zexpanded}, one has (we remark that the operator $\mathcal{L}_{\gamma}^{*}$ is real)
\begin{equation}
C_{\mathcal{L}_{\gamma}^{*}Z} = \sum_{n \in \mathbb{N}\setminus N_{0}} \mathcal{L}_{\gamma}^{*}( \sqrt{\lambda_{n}}T_{n} ) \otimes  \mathcal{L}_{\gamma}^{*}( \sqrt{\lambda_{n}}T_{n} ) = \sum_{n \in \mathbb{N}} g_{n}^{*} \otimes g_{n}^{*},
\end{equation}
from which it is immediate that
\begin{equation}
\langle C_{\mathcal{L}_{\gamma}^{*}Z} , \varphi \otimes \overline{\phi} \rangle = \sum_{n \in \mathbb{N}} ( \varphi , g_{n} )_{L^{2}(\Rd)}(\overline{\phi} , g_{n})_{L^{2}(\Rd)} = (\varphi , \phi )_{L^{2}(\Rd)}, \quad \forall\varphi , \phi \in \SchwartzRd.
\end{equation}
Hence $\mathcal{L}_{\gamma}^{*}Z$ is a White Noise.

The proof of the converse comes from a simple dimensional analysis. If we suppose that $\mathbb{N}\setminus N_{0}$ if finite and that there exists $\mathcal{L}^{*} : \mathscr{H}_{Z}'(\Rd) \to \TemperedRd$ applicable to $Z$ such that $\mathcal{L}^{*}Z$ is a White Noise, we would have from \eqref{Eq:CL*Zexpanded} the covariance equalities
\begin{equation}
C_{\mathcal{L}^{*}Z} = \sum_{n \in \mathbb{N}} e_{n}^{*} \otimes \overline{e_{n}^{*}} =  \underbrace{\sum_{n \in \mathbb{N}\setminus N_{0}}}_{\small \hbox{finite sum} \normalsize} \mathcal{L}^{*}(h_{n}) \otimes  \overline{\mathcal{L}^{*}}(h_{n}),
\end{equation}
being $(e_{n})_{n \in \mathbb{N}}$ any orthonormal basis of $L^{2}(\Rd)$. It suffices to apply these distributions to tensor products between members of the orthonormal basis $(e_{n})_{n \in \mathbb{N}}$ to see that the equalities cannot hold for an enough-large amount of these tensors. $\blacksquare$

\section{Concluding remarks}
\label{Sec:Conclusion}

We conclude by making some remarks about these results.

\begin{Rem}
\label{Rem:StrictlyPosDef}
A simple case when Theorem \ref{Theo:LZ=W} holds is when the covariance distribution of $Z$ is (strictly) positive definite, that is, $\langle C_{Z} , \varphi\otimes \overline{\varphi}\rangle = 0 $ if and only if $\varphi = 0$. In such a case $C_{Z}$ induces an inner product, and it suffice to take the completition of $\SchwartzRd$ with respect to this inner product as $\mathscr{H}_{Z}(\Rd)$, being the latter infinite dimensional. This is the case, for example, of mean-square continuous  (second-order) stationary stochastic process whose spectral measures have strictly positive density.
\end{Rem}

\begin{Rem}
\label{Rem:N=0=M}
When the GeSP $Z$ is already a mean-square continuous stochastic process, one can set $N = 0$ in \eqref{Eq:DefY}. If in addition the covariance is bounded, we can also set $M = 0$ in \eqref{Eq:CYBounded}, and the results obtained rely mainly on the Karhunen-Loève expansion of the process (with respect to the Hilbert space $\LtwoRdmu$). This is the case, for example, of mean-square continuous (second-order) stationary stochastic processes.
\end{Rem}

\begin{Rem}
\label{Rem:Gaussian}
If the real GeSP is Gaussian, all the constructed stochastic processes and random variables used in this work can be taken to be Gaussian. This comes from the stability of Gaussian processes against linear transformations. In such case,  the word ``uncorrelated'' can be replaced by ``independent'' all along this paper (for Theorem \ref{Theo:KL} it applies if the process is real).
\end{Rem}

\begin{Rem}
\label{Rem:GeneralDistributions}
The case of $Z$ being a non-tempered random distribution, that is, when $Z$ is a continuous linear operator from $\DRd \to \LtwoOmega$, being $\DRd$ the classical space of smooth functions with compact support, then the results presented here can be obtained \textit{locally}. That is, for every bounded open set $O \subset \Rd$, analogous results to Theorems \ref{Theo:Z=LW} and \ref{Theo:LZ=W}, as well as to Corollary \ref{Corol:Zexpression} can be obtained by restricting the analysis to test-functions whose supports are contained in $O$. This comes from the Regularity Theorem of General Distributions  \citep[Theorem XX1, Section 6, Chapter III]{schwartz1966theorie}, which indicates, analogously to Theorem \ref{Theo:RegularityTempered}, that every member of $\DistributionsRd$ is locally the derivative of high-enough order of a continuous function. Our arguments can be then applied to such case.
\end{Rem}

\section*{Acknowledgements}

The main of this research has been done during my stay as research and teaching assistant at the \textit{École Nationale Supérieure d'Informatique pour l'Industrie et l'Entreprise} (ENSIIE) in Évry, France, institution to which I am very grateful for its reception.

\begin{appendices}

\section{Reminders on Topological Vector Spaces and proof of Proposition \ref{Prop:ExtenstionZtoU}}
\label{Sec:AppReminders}

In this Appendix we follow mainly \citet[Chapter V]{reed1980methods} and \citet[Chapter 43]{treves1967topological}.

A (complex) Hausdorff locally convex topological vector space $E$ is a vector space endowed with a topology induced by a (supposed directed) family of  seminorms $(p_{j})_{j \in J}$ such that if $p_{j}(e) = 0$ for all $j \in J$ then $e = 0$. The induced topology consists in the weakest topology which makes continuous the sum, the scalar multiplication and the seminorms $(p_{j})_{j \in J}$. In such case, the topological dual space $E'$ consists of linear functionals $T \in E^{*}$ such that there exists $C > 0$ and $j \in J$ for which
\begin{equation}
|\langle T , e \rangle | \leq C p_{j}(e), \quad \forall e \in E.
\end{equation}
As example, the Schwartz space $\SchwartzRd$ is classically endowed with the countable family of seminorms (other equivalent seminorms can be given)
\begin{equation}
\label{Eq:SemiNormsSchwartz}
\mathcal{N}_{n}(\varphi) = \sum_{\substack{\beta \in \mathbb{N}^{d} \\ |\beta| \leq n }} \| (1+|x|^{2})^{\frac{n}{2}} D^{\beta}\varphi \|_{\infty}, \quad \forall \varphi \in \SchwartzRd, n \in \mathbb{N},
\end{equation}
where $D^{\beta}$ denotes a differential operator of order $\beta \in \mathbb{N}^{d}$ (multi-index notation), and $\| \cdot \|_{\infty}$ denotes the supremum norm.

Let $F$ be another Hausdorff locally convex topological vector space with directed family of seminorms $(q_{j'})_{j' \in J'}$. A linear operator $\mathcal{L} : E \to F $ is then continuous if and only if for every $j' \in J'$ there exists $C_{j'} > 0$ and $j \in J $ such that 
\begin{equation}
q_{j'}( \mathcal{L}(e) ) \leq C_{j'} p_{j}(e), \quad \forall e \in E.
\end{equation}
If $E \subset F $, $E$ is said to be continuously embedded in $F$ if the identity mapping $E \to E \subset F$ is continuous. This continuous embedding condition implies $F' \subset E'$. As example, $\SchwartzRd$ is continuously embedded in $L^{2}(\Rd)$, for which it is easy to prove that for some $C'>0$, $\| \varphi \|_{L^{2}(\Rd)} \leq C' \mathcal{N}_{d+1}(\varphi)$  for every $\varphi \in \SchwartzRd$.

The tensor product $E \otimes F$ can be endowed with the projective topology, which is the strongest locally convex topology which makes continuous the canonical operator $E \times F \to E \otimes F $. This topology is the one induced by the seminorms
\begin{equation}
\label{Eq:TensorSemiNorms}
( p_{j} \otimes p_{j'} )(\psi)  = \inf \lbrace   \sum_{i \in I} p_{j}(e_{i})q_{j'}(f_{i}) \ \big| \   \psi = \sum_{i \in I}e_{i}\otimes f_{i}, I \hbox{ finite}  \rbrace, \quad  \forall \psi \in E\otimes F,  (j,j') \in J\times J'.
\end{equation}
The tensor product endowed with this topology is denoted $E\otimes_{\pi}F$. In consequence, it is possible to prove that a linear functional $T \in ( E \otimes_{\pi} F )^{*}$ is in $( E \otimes_{\pi} F )'$ if and only if there exists $C > 0 $ and $(j,j') \in J \times J'$ such that
\begin{equation}
| \langle T , e \otimes f \rangle | \leq C p_{j}(e)p_{j'}(f), \quad \forall (e,f) \in E\times F.
\end{equation}
 
Now we give the proof of Proposition \ref{Prop:ExtenstionZtoU}. Let us consider $\mathscr{U}(\Rd)$ as in Proposition \ref{Prop:ExtenstionZtoU}. Since $\SchwartzRd \subset \mathscr{U}(\Rd)$ with continuous embedding, we have $\mathscr{U}'(\Rd) \subset \TemperedRd$, and it is also possible to prove using the already given definitions that
\begin{equation}
(\mathscr{U}(\Rd)\otimes_{\pi}\mathscr{U}(\Rd))' \subset (\SchwartzRd\otimes_{\pi}\SchwartzRd)' = \TemperedRdxRd,
\end{equation}
where the last equality is given by the Schwartz's Nuclear Theorem. Thus, $(\mathscr{U}(\Rd)\otimes_{\pi}\mathscr{U}(\Rd))'$ is indeed a subspace of tempered distributions over $\RdxRd$, for which it makes sense to ask if the covariance distribution $C_{Z} \in \TemperedRdxRd$ belongs to it. If it does, as it is supposed, we have hence that there exists $C''> 0 $ and some seminorm $p$ on $\mathscr{U}(\Rd)$ such that
\begin{equation}
| \langle C_{Z} , \varphi \otimes \overline{\phi} \rangle | \leq C'' p(\varphi)p(\phi), \quad \forall \varphi,\phi \in \SchwartzRd.
\end{equation}
This implies, by definition of the covariance distribution, that
\begin{equation}
\label{Eq:ZcontinuousU}
 \mathbb{E}(|Z(\varphi)|^{2}) \leq C''p(\varphi)^{2}, \quad \forall \varphi \in \SchwartzRd, 
\end{equation}
hence $Z$ is continuous as a linear operator from $\SchwartzRd$ to $\LtwoOmega$ with $\SchwartzRd$ being endowed with the subspace topology induced by $\mathscr{U}(\Rd)$. Proposition \ref{Prop:ExtenstionZtoU} follows hence from a typical extension of continuous linear functions to the completition of the domain, which we will make precise for sake for completeness (\citet[Theorem 5.1]{treves1967topological}, considering only the linear and continuous case). We recall that we have supposed that $\SchwartzRd$ is dense in $\mathscr{U}(\Rd)$.

\begin{Theo}
Let $E,F$ be two Hausdorff topological vector spaces, $D$ a dense subspace of $E$, and let $f : D \to F$ be a continuous linear mapping from $D$ to $F$. Then, if $F$ is complete, there exists a unique continuous linear mapping $\hat{f} : E \to F $ which extends $f$, that is, such that $\hat{f}(x) = f(x)$ for all $x \in D$.
\end{Theo} 
Proposition \ref{Prop:ExtenstionZtoU} follows then simply using $f = Z$, $D = \SchwartzRd$, $E = \mathscr{U}(\Rd)$, and $F = \LtwoOmega$, the latter being complete since it is a Hilbert space. The definition of $\mathcal{L}^{*}Z$ as a GeSP is immediately concluded from the continuous extension of $Z$ to $\mathscr{U}(\Rd)$. 
\end{appendices}

\bibliography{mibib}
\bibliographystyle{apacite}

\end{document}